\documentclass{article}
\usepackage{graphics}
\usepackage{amsthm}
\usepackage{amssymb}
\usepackage{subfigure}
\usepackage[fleqn]{amsmath}
\usepackage{float}
\usepackage{algorithm}
\usepackage{algorithmic}
\usepackage{amsfonts,mathrsfs,amsmath,amssymb,amsthm}
\usepackage{graphicx,hyperref,lipsum,psfrag,stfloats,url}
\usepackage{array,footnote,longtable,subfigure}
\usepackage{multirow,multicol,verbatim}
\usepackage[usenames]{color}
\usepackage{color}
\usepackage{epstopdf}

\newtheorem{theorem}{Theorem}
\newtheorem{corollary}{Corollary}
\newtheorem{definition}{Definition}

\newtheorem{lemma}{Lemma}

\begin{document}
\title{State Diagrams  of a Class of Singular LFSR and Their Applications to the Construction of de Bruijn Cycles}
\author{XiaoFang Wang, YuJuan Sun and WeiGuo Zhang
\\ ISN Laboratory, Xidian University, Xi'an 710071, China\\
e-mail: xiaofangwang77@163.com}
\date{}
\maketitle

\begin{abstract}
The state diagrams  of a class of singular linear feedback shift registers (LFSR)  are discussed.
It is shown that the state diagrams of the given LFSR have  special structures.
An algorithm is presented to construct a new class of de Bruijn cycles from the state diagrams of these singular LFSR.
\end{abstract}


\section{Introduction}
\label{sec1-1}
Pseudo-random sequences that come from feedback shift registers (FSR) are basic components for stream ciphers in cryptography, and  FSR are also important components in modern communications due to their ease and efficiency in implementation.
The problems on FSR have been important challenges in mathematics and they were extensively studied (see \cite{Golomb67,GolombGong05,Alhakim2012Spans}, and the references therein).

As a kind of special FSR sequences, de Bruijn sequences have attracted much attentions due to their good randomness properties \cite{Fredricksen1978Necklaces,Hauge1996De,Hauge1996On,Alhakim2012Spans,Mayhew2002Extreme,MandalGong131,ChangEzerman161,ChangEzerman162,Sawada2016A}.
A comprehensive survey of previous works on this subject can be found in \cite{Fredricksen82}.

The sequences generated by nonsingular linear feedback shift registers (LFSR)  have been analyzed  using some methods from finite fields \cite{Golomb67,GolombGong05}.
Their properties can be explained by the description of their geometric structures.
For examples, maximum-length LFSR, pure cycling registers and pure summing registers have special cycle structures \cite{EtzionLempel84}.
Nonsingular LFSR can be used to generate de Bruijn cycles by the cycle-joining method \cite{Golomb67}.
The applications of this method requires the full knowledge of the cycle structures and adjacency graphs of the original LFSR
\cite{HemmatiSchilling84,LiHelleseth141,LiHelleseth142,LiHelleseth16,LiLin16}.
However, the cycle-joining method does not work  to construct de Bruijn cycles from singular LFSR since the state diagrams don't only contains disjoint cycles.
The number of  singular LFSR is the same as that of nonsingular LFSR \cite{Golomb67}.
But unfortunately, for a fixed singular LFSR, the geometry structure of its state diagram is far less developed, let alone its other properties.

The purpose of this paper is to characterize the state diagrams of a class of singular LFSR  and to construct some new de Bruijn sequences.
The remainder of this paper is organized as follows.
In Section~\ref{sec1-2}, we introduce some basic conceptions and related results.
In Section~\ref{sec1-3}, the state diagrams of a class of singular LFSR are completely determined, and the adjacency graphs of the LFSR are also given.
In Section~\ref{sec1-4}, we construct a new family of de Bruijn cycles.
In Section~\ref{sec1-5}, we provide an example to illuminate the algorithms given in the previous sections.
Section \ref{sec1-6} concludes the paper.

\section{Preliminaries}
\label{sec1-2}
For a positive integer $n$, let $GF(2)^n$ be the $n$-dimensional vector space over $GF(2)$, where $GF(2)$ is the finite field with two elements.

\subsection{FSR Sequences}
\label{subsec1-1}
An \emph{$n$-stage FSR} is a circuit arrangement consisting of $n$ binary storage elements (called stages) regulated by a clock.
The stages are labeled from $1$ to $n$.
A \emph{state of an FSR} is an $n$-tuples  vector  $(x_1, x_2, \ldots , x_{n})$, where $x_i$  indicates  the content of stage $i$.
At each clock pulse, for some integer $i \geq 1$, the state $A_i = (a_i , a_{i+1}, \ldots , a_{i+n-1})$ is updated to
\begin{equation}
\label{equ1-1}
  A_{i+1 }= (a_{i+1}, \ldots , a_{i+n-1}, F(a_i , a_{i+1}, \ldots , a_{i+n-1})),
\end{equation}
where $F: GF(2)^n\mapsto GF(2)$ is the feedback function of the FSR.
The initially loaded content $A_1 = (a_1, a_2, \ldots , a_{n})$ is called an \emph{initial state}.
The feedback function $F$ can induce a \emph{state operation} $\Gamma: GF(2)^n\mapsto GF(2)^n$ with $\Gamma(A_{i})=A_{i+1}$.
After consecutive clock pulses, the FSR outputs a sequence $a = (a_1, a_2, a_3, \ldots )$ called an output sequence of the FSR.
Therefore, the sequence $a = (a_1, a_2, a_3, \ldots )$  satisfies the recursive relationship
\begin{equation}
\label{equ1-2}
a_{i+n} =F(a_i , a_{i+1}, \ldots , a_{i+n-1})
\end{equation}
for all $i \geq 1$ and a given initial state $A_1 = (a_1, a_2, \ldots , a_{n})$.
We denote by $\Omega( F )$ the set of all $2^n$ sequences generated by $F(x_1, x_2, \ldots , x_{n})$.

An FSR is called an \emph{LFSR} if its feedback function $F(x_1, x_2, \ldots , x_{n})$ is linear.
The output sequences of an LFSR are called LFSR sequences.
A \emph{de Bruijn sequence of order $n$} is a binary sequence $a = (a_1, a_2, a_3,\ldots)$  of period $2^n$ such that all $n$-tuples $(a_i, a_{i+1}, \ldots,a_{i+n-1}), i=1,2,\ldots,2^n$, are pairwise distinct, where $n$ is a positive integer.

Let $\mathbb{V}(GF(2))$ be the set of all binary infinite sequences.
For $a =(a_1, a_2, a_3, \ldots) \in \mathbb{V}(GF(2))$, we define the \emph{shift operator} $L$ on $\mathbb{V}(GF(2))$ by
$$La =L(a_1, a_2, a_3, \ldots)  = (a_2, a_3, a_4, \ldots) \in \mathbb{V}(GF(2)).$$
Generally, we have $L^ma  = (a_{m+1}, a_{m+2},a_{m+3}, \ldots)$ for any integer $m$ with $m \geq 1$, and denote $L^0a = a$.
A sequence $a =(a_1, a_2, a_3, \ldots) \in \mathbb{V}(GF(2))$ is called periodic if there exists a positive integer $r$ such that $L^r a = a$,  i.e.,   $a_{i+r} = a_i$ for each integer $i \geq 1$.
The smallest positive integer $r_0$ satisfying this property is called the period of $a$, denoted by $per(a)$.
Specially, we can define the shift operator $L$ on $GF(2)^n$.
For an $n$-stage FSR, let $A_1 =(a_1, a_2, \ldots,a_{n}) \in GF(2)^n$, we have
\begin{equation}
\label{equ1-3}
  L^iA_1=A_{i+1}=(a_{i+1}, a_{i+2},\ldots , a_{i+n}),~~ where~~ i\geq 0.
\end{equation}

An $n$-stage FSR and its feedback function $F(x_1, x_2, \ldots , x_{n})$ are said to be \emph{nonsingular} if the state operation $\Gamma$ is one-to-one, i.e., for any two states $A$ and $B$, $\Gamma A = \Gamma B$ implies $A = B$.
It is well known that all sequences in $\Omega( F )$ are periodic if and only if the feedback function $F$ is nonsingular, i.e., $F$ can be written as $F(x_1 , x_2, \ldots , x_{n})=x_1+F_1(x_{2}, x_{3},\ldots , x_{n})$, where $F_1$ is a function from $GF(2)^{n-1}$ to  $GF(2)$.
If the state operation $\Gamma$ is not one-to-one, then the FSR and its feedback function $F$ are said to be \emph{singular}, i.e., at least there exist two different states $A^{\prime}$ and $B^{\prime}$ satisfying  $\Gamma A ^{\prime}= \Gamma B^{\prime}$.
When $F$ is singular, i.e., $F$ can not be written as $F(x_1 , x_{2}, \ldots , x_{n})=x_1+F_1(x_{2}, x_{3},\ldots , x_{n})$, and all the sequences in $\Omega( F )$ are ultimately periodic.
The feedback function $F(x_{1},x_{2},\ldots, x_{n})=x_{n-1}+x_n \ (n\geq3)$ discussed in this paper is obviously singular.

\subsection{State diagrams}
\label{subsec1-2}
A directed graph  $D$ is an ordered pair $(V(D),A(D))$ consisting of a set $V(D)$ of vertices and a set $A(D)$  of arcs together with an incidence function $\psi_D$, such  that each arc of $A(D)$ is an ordered pair of vertices of $V(D)$.
If $\mathfrak{a}~:u\rightarrow v$ is an arc, i.e., $\psi_D(\mathfrak{a}) = (u, v)$, then $\mathfrak{a}$ is said to join $u$ to $v$, and we also say that $u$ dominates $v$.
The vertex $u$ is the tail of $\mathfrak{a}$, and the vertex $v$ is its head.
They are the two ends of $\mathfrak{a}$.

For an $n$-stage FSR with the feedback function $F(x_1, x_2, \ldots , x_{n})$, its \emph{state diagram} $G_F$ is a directed graph with $2^n$ vertices, each vertex labeled with a unique binary $n$-length vector.
For two vertices $Y^n =(y_1, y_{2}, \ldots, y_{n})$ and $Z^n =(z_1, \ldots, z_{n-1}, z_{n})$ in $G_F$, an arc is drawn from $Y^n$ to $Z^n$ if and only if $( y_2, y_3, \ldots, y_{n})=(z_1, z_2,  \ldots, z_{n-1})$ and $z_n =F(y_1, y_2,  \ldots, y_{n})$.
In other words, an $n$-stage FSR state diagram $G_F$ is a directed graph with $2^n$ vertices such that there is an  arc $Y^n\rightarrow Z^n$ if and only if $L(Y^n) = Z^n$, where $L$ is the shift operator.
In this case, we say that $Y^n$ leads to $Z^n$, $Y^n$ is a \emph{predecessor} of $Z^n$ and $Z^n$ is the \emph{successor} of $Y^n$.
We also say that the vertex $Y^n$ is the tail of  the arc $Y^n\rightarrow Z^n$, and the vertex $Z^n$ is its head.
In particular, if $Y^n$ is a predecessor of $Y^n$,  then the state diagram has a \emph{loop} at $Y^n$.
For any given $n$-length vector $Y^n$, there are two possible predecessors for $Y^n$, and $Y^n$ also has two possible successors.
Given an $n$-stage FSR, each state has  only one successor.

A \emph{cycle} $C = (Y^n_{1}, Y^n_{2}, \ldots ,Y^n_{k})$ of length $k$ in the state diagram of an $n$-stage FSR  is a cyclic sequence of $k$ distinct states $Y^n_{1}, Y^n_{2}, \ldots ,Y^n_{k}$ such that $Y^n_{t}$ leads to $Y^n_{t+1}$ for each $1\leq t \leq k-1$ and $Y^n_{k}$ leads to $Y^n_{1}$, where $k\geq1$.
A $k$-length cycle is denoted also by $C(k)$.
A convenient representation of a cycle of length $k$ is  a ring sequence $[y_1,y_2,\ldots,y_{k}]$, where $y_i$ is the first component of $Y^n_{i}$ for $1\leq i\leq k$.
From the definition of a cycle, all the cyclic shifts of $[y_1,y_2,\ldots,y_{k}]$ denote the same cycle.

For a fixed singular FSR, the state diagram $G_F$ can be denoted by $G_F=G_1\cup G_2\cup\ldots\cup G_t$, where $G_i$, $i=1,2,\ldots,t$, is a connected component.

The maximum length of cycles in the state diagram of an $n$-stage FSR is $2^n$.
In this case, the FSR is said to be a maximum-length FSR and the cycle is  an $n$-order de Bruijn cycle (full cycle).

For a state $Y^n =(y_1,y_2,\ldots, y_{n})$ of an $n$-stage FSR, its \emph{conjugate} $\widehat{Y^n}$ and \emph{companion} $\widetilde{Y^n}$ are defined as
$$\widehat{Y^n}=(y_1+1, y_2,  \ldots, y_{n})$$
and
$$\widetilde{Y^n}=(y_1, y_2,  \ldots, y_{n}+1).$$
Two states $Y^n$ and $\widehat{Y^n}$ form a conjugate pair, and a conjugate pair is denoted  by $conj(Y^n,\widehat{Y^n})$.
Two states $Y^n$ and $\widetilde{Y^n}$ constitute a companion pair, and a conjugate pair is denoted  by $comp(Y^n,\widetilde{Y^n})$.
Given an $n$-length vector $Y^n =(y_1, y_2, \ldots, y_{n})$, its two possible predecessors form a conjugate pair, and its two possible successors constitute a companion pair.

For a nonsingular FSR, the problem of determining conjugate pairs among cycles in a state diagram $G_F$ leads to the definition of  adjacency graph.
For a fixed nonsingular FSR   with feedback function $F(x_1, x_2, \ldots , x_{n})$, its   adjacency graph is an undirected weighted graph where the vertices correspond to the cycles $C_1,C_2,\ldots,C_s$ in $G_F$.
Two cycles $C_i$ and $C_j$ are \emph{adjacent} if they are disjoint and there exists a state $Y^n$ on cycle $C_i$  whose conjugate $\widehat{Y^n}$ is on cycle $C_j$, where $1\leq i,j\leq s$.
Then we also say that cycle $C_j$ is the adjacent cycle of  cycle $C_i$.
The classic idea of the cycle-joining method is that two adjacent cycles $C_i$ and $C_i$ are joined into a single cycle when the successors of $Y^n$ and $\widehat{Y^n}$ are interchanged.
When two cycles can be connected by exactly $m$ edges, it is convenient to denote the $m$ edges by an edge labeled with an integer $m$, where $m\geq 2$.
The weighted graph obtained in this way is called the  adjacency graph of the nonsingular  FSR  with the feedback function $F(x_1, x_2, \ldots , x_{n})$ and denoted by $\mathcal {G}(F)$.

For any singular FSR, we can obtain its adjacency graph $\mathcal {G}(F)$ by making the  connected components $G_1,G_2,\ldots,G_t$ corresponding to vertices similarly.

\subsection{Perfect binary directed trees}
\label{subsec1-3}
Each acyclic connected graph is a tree.
The top vertex in a tree is its root.
A  leaf (or terminal vertex) is a vertex with only one adjacent vertex.
The depth of a vertex is the number of edges from this vertex to the tree's root vertex.
The binary tree is a connected acyclic graph, and   each vertex has at most two adjacent vertices.
If a directed graph is an acyclic connected graph when  the direction is ignored, then we call this directed graph as a directed tree.
If the tree is a binary tree when there is no consideration to the direction,  then we call this directed graph as a binary directed tree.

\begin{definition}\label{def1-1}
A perfect binary directed tree is a binary directed tree in which all interior vertices have two adjacent vertices and all leaves have the same depth or same level.
The depth of a perfect  binary directed tree is the length of path from one leaf  to the root.
A perfect  binary  directed tree with depth $k$ is denoted by $T_k$.
Obviously, there are $(2^{k+1}-1)$ vertices in $T_k$.
\end{definition}

Given a  digraph $D$, denote  all vertices (or arcs) of $D$ as $V(D)$ (or $A(D)$).
If  $v\in V(D)$ is a head or a tail of  $a\in A(D)$, then we call $a$ to be incident with $v$.
Let all arcs incident with  $v$ be denoted by $A(v)$, then $D\setminus v$ can be induced from $D$ such that $V(D\setminus v)=V(D)\setminus\{v\}$ and $A(D\setminus v)=A(D)\setminus A(v)$.
If $C = (v_1,v_2,\ldots,v_k)$ is a $k$-length cycle and $V(C)=\{v_1,v_2,\ldots,v_k\}$, $D\setminus C$ can be induced from $D$ such that $V(D\setminus C)=V(D)\setminus\{v_1,v_2,\ldots,v_k\}$ and $A(D\setminus C)=A(D)\setminus(A(v_1)\cup A(v_2)\cup \ldots \cup A(v_k))$.

\section{State diagrams  of LFSR  with the feedback function $F(x_{1},x_{2},\ldots, x_{n})=x_{n-1}+x_n $}
\label{sec1-3}
We next consider the singular $n$-stage LFSR with the feedback function of the form
\begin{equation}
\label{1-4}
 F(x_{1},x_{2},\ldots, x_{n})=x_{n-1}+x_n,~(n\geq3).
\end{equation}

\begin{lemma}
\cite{Golomb67}
\label{state-diagram1-1}
 Given an  FSR  with feedback function $F(x_{1},x_{2},\ldots, x_{n})$, the state diagram $G_F$ only has finite connected components $G_1,G_2,\ldots,G_t$.
 Every connected component only  contains  unique cycle.
 In other words, every connected component is a cycle or a cycle with some branches.
\end{lemma}

\begin{lemma}
\cite{Golomb67}
\label{state-diagram1-2}
 Let $F(x_{1},x_{2},\ldots, x_{n})$ be the feedback function of a fixed $n$-stage FSR, and  $a= (a_1, a_2, a_3, \ldots )\in \Omega( F )$ be  the output sequence of a given initial state $A_1^n = (a_1, a_2, \ldots , a_{n})$.
 Then  $A_1^n$ belongs to unique connected component $G_i$ in the state diagram $G_F$.
 Furthermore, if we delete some initial terms of $a$,  a periodic sequence whose  period is the length of the cycle $C_i$ can be obtained, where $C_i$ is the unique cycle of $G_i$.
\end{lemma}

For $(a,b,c)\in GF(2)^3$, let
$$S_{abc}^n=\{Y^n=(y_1,y_2, \ldots, y_{n})\mid (y_{n-2}, y_{n-1},  y_{n})=(a,b,c)\}.$$
Let $n$-length vectors of zeros and ones be denoted by $\mathbf{0}^n=(0,0,\ldots,0)$ and $\mathbf{1}^n=(1,1,\ldots,1)$.

\begin{theorem}
\label{state-diagram1-3}
 Given an  $n$-stage singular LFSR with feedback function $ F(x_{1},x_{2},\ldots, x_{n})=x_{n-1}+x_n$,  the state diagram  $G_F$ only has two connected components $G_1$ and $G_2$, where  $G_1$ contains  the loop $[0]$ and  $G_2$  contains  a $3$-length cycle $[0,1,1]$.
\begin{proof} For any state $Y^n=(y_1,y_2,\ldots,y_n)$  of
\begin{equation*}
S_{100}^n=\{Y^n=(y_1,y_2, \ldots,y_{n-2}, y_{n-1},  y_{n}) \mid (y_{n-2}, y_{n-1},  y_{n})=(1,0,0)\},
\end{equation*}
we have
\begin{eqnarray*}
   LY^n &=& (a_2,~\ldots,~a_{n-3},~1,0,0,0) \\
   L^2Y^n &=& (a_3,\ldots,a_{n-3}, 1, 0, 0, 0, 0) \\
      &\ldots&  \\
   L^{n-2}Y^n &=& (0,\ldots,0, 0, 0, 0, 0, 0)= \mathbf{0}^n.
\end{eqnarray*}
Similarly, for any  state $Y^n=(y_1,y_2,\ldots,y_n)$ belonging to  $S_{000}^n$,  we have
\begin{equation*}
 LY^n=(y_2,y_3,\ldots,0,0,0,0),~\ldots,~L^iY^n=\mathbf{0}^n,~0\leq i\leq n-3.
\end{equation*}
Obviously, when $\mathbf{0}^n$ is the initial state, its successor is $\mathbf{0}^n$ itself.
That is, there is a loop at  $\mathbf{0}^n$.
Then every  state $Y^n \in S_{100}^n\bigcup S_{000}^n $ always leads to $\mathbf{0}^n$ and  belongs to the same connected component $G_1$, and this connected component contains the loop $[0]$.

Let $v_1=(0,1,1,0,1,1,\ldots)$, $v_2=(1,1,0,1,1,0,\ldots)$ and $v_3=(1,0,1,1,0,1,\ldots)\in GF(2)^n$.
Since the feedback function is $F(x_{1},x_{2},\ldots, x_{n})=x_{n-1}+x_n $, we can get the following relations through calculation: $Lv_1= v_2$, $Lv_2= v_3$ and $Lv_3= v_1$.
In other words, $v_1\rightarrow v_2 \rightarrow v_3\rightarrow v_1$, i.e., $[0,1,1]$ is a $3$-length cycle.
It is not difficult to verify that for any state $Y^n\not\in S_{100}^n\bigcup S_{000}^n$, $Y^n$ will ultimately go into the cycle $[0,1,1]$, which implies $Y^n$ belongs to the other connected component $G_2$ and $G_2$ contains  the $3$-length cycle $[0,1,1]$.
\end{proof}
\end{theorem}

\begin{theorem}
\label{state-diagram1-4}
With the same notations as above.
$G_1\setminus[0]$  is a perfect  binary directed tree  $T_{n-3}^{(1)}$, and $G_2\setminus[0,1,1]$  contains three perfect  binary directed trees  $T_{n-3}^{(2)}$, $T_{n-3}^{(3)}$ and $T_{n-3}^{(4)}$.
\begin{proof}
Note that the $n$-stage LFSR is singular, we have that a pair conjugate states have the same successor.
By Lemma~\ref{state-diagram1-1} and Theorem~\ref{state-diagram1-3},  $G_1\setminus[0]$ and $G_2\setminus[0,1,1]$ are four binary directed trees.
For each $Y^n\in S_{100}$, we have $L^{n-3}Y^n=\widehat{\mathbf{0}^n}$.
Since the length of the directed path from any leaf $Y^n\in S_{100}$ to $\widehat{\mathbf{0}^n}$ is $n-3$ ($\widehat{\mathbf{0}^n}$ is the root of this tree),  this  binary directed tree $G_1\setminus[0]$ is a perfect  binary directed tree $T_{n-3}^{(1)}$.
Similarly, $G_2\setminus[0,1,1]$  contains three perfect  binary directed trees $T_{n-3}^{(2)}$, $T_{n-3}^{(3)}$ and $T_{n-3}^{(4)}$, whose roots are $\widehat{v_1}$, $\widehat{v_2}$ and $\widehat{v_3}$, respectively.
\end{proof}
\end{theorem}

\begin{corollary}
\label{state-diagram1-5}
Let $F(x_{1},x_{2},\ldots, x_{n})=x_{n-1}+x_n \ (n\geq3)$ be the feedback function of an $n$-stage singular LFSR.
Then the adjacency graph  $\mathcal {G}(F)$ only has two isolated vertices $G_1$ and $G_2$.
\end{corollary}

\section{A new class of de Bruijn cycles}
\label{sec1-4}

In this section, we introduce a new method to obtain a new class of de Bruijn cycles by modifying the state diagrams of LFSR with the feedback function $F(x_{1},x_{2},\ldots, x_{n})=x_{n-1}+x_n \ (n\geq3)$.
Note that by Corollary ~\ref{state-diagram1-5} the  cycle-joining method does not work in this case.

Given an  $n$-stage singular  LFSR  with the feedback function $F(x_{1},x_{2},\ldots, x_{n})=x_{n-1}+x_n \ (n\geq3)$, we next present an algorithm to construct  the cycles without branches from its state diagram $G_F$.

Firstly,  all states on the state diagram $G_F$ are  classified.
The set consisting of all states on  the two cycles  $[0]$  and $[0,1,1]$ is denoted as $S_0^n$, i.e., $S_0^n=\{(0,1,1,0,1,1,\ldots),$ $ (1,1,0,1,1,0,\ldots),(1,0,1,1,0,1,\ldots),\mathbf{0}^n\}$.
The set of all leaves is denoted as $S^n$ with $S^n=S_{100}^n\cup S_{001}^n\cup S_{010}^n\cup S_{111}^n$.
An $l$-length  directed path $[Y_{s,t}^n\rightarrow LY_{s,t}^n\rightarrow\ldots\rightarrow L^{l}Y_{s,t}^n]$ is written as  $p_{s,t}$.
For a directed  path $p_{s,t}$,  the tail of $p_{s,t}$ is denoted by $T(p_{s,t})=Y_{s,t}^n$, then  the set of states remaining on the path except the tail $Y_{s,t}^n$ is denoted as $U_{s,t}$, i.e., $U_{s,t}=\{LY_{s,t}^n,\ldots,L^{l}Y_{s,t}^n\}$.
In the state diagram $G_F$, if a state has two predecessors and one successor, then we call this state a trigeminal vertex.
Obviously, $S_0^n$ and $U_{s,t}$ are all sets of trigeminal vertices.

\begin{algorithm}[h!]
\caption{Generation of the two sets $P$ and $C$ determined by the feedback function $F(x_{1},x_{2},\ldots, x_{n})=x_{n-1}+x_n \ (n\geq3)$}
\label{alg1-1}
\begin{algorithmic}[1]
\REQUIRE The feedback function $F(x_{1},x_{2},\ldots, x_{n})=x_{n-1}+x_n $, and the  state $Y_{i,j}^n\in S^n=S_{100}^n\cup S_{001}^n\cup S_{010}^n\cup S_{111}^n$
  \ENSURE The  set of paths $P$ and the set of cycles $C$
  \STATE Set $i=1$, $S_0^n=\{(0,1,1,0,1,1,\ldots), (1,1,0,1,1,0,\ldots), (1,0,1,1,0,1,\ldots), \mathbf{0}^n\}$, $T(P_{0})=\emptyset$, $U_{0,0}=\emptyset$
  \WHILE{$S^n\setminus\bigcup_{0\leq s\leq i-1}T(P_{s})\neq\emptyset$}
   \STATE $j=1$
  \WHILE{$Y_{i,j}^n\neq Y_{i,1}^n$ or  $j=1$}
  \IF{$Y_{i,j}^n$ be the firstly-appearing state of   one of the four set $S_{100}^n$, $S_{001}^n$, $S_{010}^n$ and $S_{111}^n$}
  \STATE $l_{i,j}=n-2$
  \STATE $p_{i,j} \gets [Y_{i,j}^n\rightarrow LY_{i,j}^n\rightarrow\ldots\rightarrow L^{l_{i,j}}Y_{i,j}^n]$
  \STATE $Y_{i,j+1}^n \gets \widetilde{L^{l_{i,j}+1}Y_{i,j}^n}$
  \ELSE
  \STATE $l_{i,j}$ be the least integer such that $L^{l_{i,j}}Y_{i,j}^n\in \bigcup_{k_1=0}^{i-1}\bigcup_{k_1=0}^{j-1}U_{k_1,k_2}$
  \IF  {$l_{i,j}>1$}
  \STATE $p_{i,j} \gets [Y_{i,j}^n\rightarrow LY_{i,j}^n\rightarrow\ldots\rightarrow L^{l_{i,j}-1}Y_{i,j}^n]$
  \STATE $Y_{i,j+1}^n \gets \widetilde{L^{l_{i,j}}Y_{i,j}^n}$
   \ELSE
   \STATE  $p_{i,j} \gets [Y_{i,j}^n]$, i.e., the single state $Y_{i,j}^n$ is considered as a path of length $0$
   \STATE $Y_{i,j+1}^n \gets \widetilde{L^{l_{i,j}}Y_{i,j}^n}$
   \ENDIF
  \ENDIF
  \STATE  $j \gets j+1$
  \ENDWHILE
  \STATE  $P_i=\{p_{i,1},p_{i,2},\ldots,p_{i,j-1}\}$
  \STATE  $C_i: [p_{i,1}\rightarrow p_{i,2}\rightarrow\cdots\rightarrow p_{i,j-1}\rightarrow p_{i,1}]$
  \STATE   $T(P_{i})=\{T(p_{i,k})|1\leq k\leq j-1, T(p_{i,k})=Y_{i,k}^n\}$
  \STATE  $i \gets i+1$
  \ENDWHILE
\RETURN {$P=\bigcup_{1\leq t\leq i-1}P_t$, $C=\{C_1,C_2,\ldots,C_{i-1}\}$}
\end{algorithmic}
\end{algorithm}

Choose a  state $Y_{1,1}^n$ from $S^n$ randomly. If  $l_{1,1}$ is the least integer such that $L^{l_{1,1}}Y_{1,1}^n\in S_0^n$, then the directed path  $[Y_{1,1}^n\rightarrow LY_{1,1}^n\rightarrow\ldots\rightarrow L^{l_{1,1}}Y_{1,1}^n]$ is denoted by $p_{1,1}$.

Based on $Y_{1,1}^n$, $l_{1,1}$ and $p_{1,1}$, select $Y_{1,2}^n$, and search for $l_{1,2}$ and $p_{1,2}$.
In general, based on $Y_{1,j}^n$, $l_{1,j}$ and $p_{1,j}$, select $Y_{1,j+1}^n$ using the  following method, and search for $l_{1,j+1}$ and $p_{1,j+1}$.
If $L^{l_{1,j}}Y_{1,j}^n\in S_0^n$, then set  $Y_{1,j+1}^n=\widetilde{L^{l_{1,j}+1}Y_{1,j}^n}$; If $L^{l_{1,j}}Y_{1,j}^n\in \bigcup_{k=0}^{j-1}U_{1,k}$, then set $Y_{1,j+1}^n=\widetilde{L^{l_{1,j}}Y_{1,j}^n}$.
Let  $l_{1,j+1}$ be  the least integer  such that $L^{l_{1,j+1}}Y_{1,j+1}^n\in S_0^n\bigcup_{k=0}^{j}U_{1,k}$. If $L^{l_{1,j+1}}Y_{1,j+1}^n\in S_0^n$, then the directed path  $[Y_{1,j+1}^n\rightarrow LY_{1,j+1}^n\rightarrow\ldots\rightarrow L^{l_{1,j+1}}Y_{1,j+1}^n]$ is denoted by $p_{1,j+1}$; otherwise,  $L^{l_{1,j+1}}Y_{1,j+1}^n\in \bigcup_{k=0}^{j}U_{1,k}$, then the directed path $[Y_{1,j+1}^n \rightarrow LY_{1,j+1}^n \rightarrow\ldots\rightarrow L^{l_{1,j+1}-1}Y_{1,j+1}^n]$ is denoted by $p_{1,j+1}$.
If $l_{1,j}=1$, then the single state $Y_{1,j}^n$ is considered as a path of length $0$, i.e., $p_{1,j}: [Y_{1,j}^n]$.
Repeat the directed path searching in the state diagram until the state $Y_{1,s_1+1}^n$ returns to $Y_{1,1}^n$.
Set $P_1=\{p_{1,1},p_{1,2},\ldots,p_{1,s_1}\}$ and $T(P_{1})=\{T(p_{1,j})|1\leq j\leq s_1\}$.
Connect the  head of $p_{1,v}$  with the  tail  of $p_{1,v+1}$ in turn for $1\leq v\leq s_1-1$, and  connect the  head of $p_{1,s_1}$  with the  tail  of $p_{1,1}$, then the  directed cycle $C_1$ is obtained, i.e. $C_1: [p_{1,1}\rightarrow p_{1,2}\rightarrow\ldots\rightarrow p_{1,s_1}\rightarrow p_{1,1}]$.

Based on the set of directed paths $P_1$ and the cycle $C_1$, the new directed path sets $P_2$ and the new cycle $C_2$ can be  searched.
In general, based on $P_i$ and  $C_i$, search for $P_{i+1}$ and  $C_{i+1}$ in the following way.
Choose an $n$-stage state $Y_{i+1,1}^n$ from $S^n\setminus \bigcup_{1\leq m\leq i}T(P_{m})$ randomly.
Let   $l_{i+1,1}$ be  the least integer  such that $L^{l_{i+1,1}}Y_{i+1,1}^n\in S_0^n\bigcup_{m\leq i}U_{m,k}$.
If $L^{l_{i+1,1}}Y_{i+1,1}^n\in S_0^n$, then the directed path  $[Y_{i+1,1}^n\rightarrow LY_{i+1,1}^n\rightarrow\ldots\rightarrow L^{l_{i+1,1}}Y_{i+1,1}^n]$ is denoted by $p_{i+1,1}$; If $L^{l_{i+1,1}}Y_{i+1,1}^n\in \bigcup_{m\leq i}U_{m,k}$, then the directed path $[Y_{i+1,1}^n \rightarrow LY_{i+1,1}^n \rightarrow\ldots\rightarrow L^{l_{i+1,1}-1}Y_{i+1,1}^n]$ is denoted by $p_{i+1,1}$.
Based on $Y_{i+1,1}^n$, $l_{i+1,1}$ and $p_{i+1,1}$, select $Y_{i+1,2}^n$ recursively, and search for $l_{i+1,2}$ and $p_{i+1,2}$.
Repeat the directed path searching in the state diagram until the state $Y_{i+1,s_{i+1}+1}^n$ returns to $Y_{i+1,1}^n$.
If $l_{i+1,j}=1$, then the single state $Y_{i+1,j}^n$ is considered as a path of length $0$, i.e., $p_{i+1,j}: [Y_{i+1,j}]^n$.
Set $P_{i+1}=\{p_{i+1,1},p_{i+1,2},\ldots,p_{i+1,s_{i+1}}\}$ and $T(P_{i+1})=\{T(p_{i+1,j})|1\leq j\leq s_{i+1}\}$.
Connect the  head of $p_{i+1,v}$  with the  tail  of $p_{i+1,v+1}$ in turn for $1\leq v\leq s_{i+1}-1$, and  connect the  head of $p_{i+1,s_{i+1}}$  with the  tail  of $p_{i+1,1}$. Then the  directed cycle $C_{i+1}$ is obtained, i.e. $C_{i+1}: [p_{i+1,1}\rightarrow p_{i+1,2}\rightarrow\ldots\rightarrow p_{i+1,s_{i+1}}\rightarrow p_{i+1,1}]$.

Repeat the directed cycle searching in the state diagram until the state set
$$S^n\setminus\bigcup_{1\leq m\leq t}T(P_{m})=\emptyset.$$
Then output the  set of paths $P=\bigcup_{1\leq m\leq t}P_m$ and the set of cycles $C=\{C_1,C_2,\ldots,C_t\}$.

\begin{theorem}
\label{state-diagram1-6}
Given an  LFSR  with feedback function $F(x_{1},x_{2},\ldots, x_{n})=x_{n-1}+x_n \ (n\geq3)$,  we repeat the directed path searching in the state diagram as in Algorithm \ref{alg1-1}.
Then  the original state diagram turns into a branchless one.
\begin{proof}
Noticing the  $n$-stage LFSR  with the feedback function $F(x_{1},x_{2},\ldots, x_{n})=x_{n-1}+x_n $, all leaves and all trigeminal vertices are completely determined.

Note  that each directed path found by Algorithm \ref{alg1-1} is a directed path with a leaf as its tail and a trigeminal vertex as its head.
We also notice that any two paths found by Algorithm \ref{alg1-1} are different in both tails and heads.
The process of finding  directed paths and connecting them as cycles in Algorithm \ref{alg1-1}, in fact, modifies the successor of one of the two predecessor states of a given trigeminal vertex to its companion state.
We also note that in this state diagram, any pair of companion states must have one  as a leaf and the other as a trigeminal vertex.
For the  head of each directed path, you can modify its successor to the companion state of the original successor, that is, the next path can be found such that its tail  is connected to  the head of the previous path.
At the same time, the tail of each directed path can be used as the modified successor of one of the two predecessor states of its companion state, i.e., there must be a path  such that its head  is connected to  the tail of this path.
So, all the leaves and trigeminal vertices in the original state diagram have disappeared.
Then the original state diagram turns into a branchless one.
\end{proof}
\end{theorem}

Algorithm \ref{alg1-1} gives the steps to generate the required  two sets $P$ and $C$.
The implementation of Algorithm \ref{alg1-1} requires the adjacency relation of state diagram.
When the feedback function is given, each state has a unique successor.
The time complexity of this algorithm is $O(2^{n+1})$.

\begin{lemma}
\cite{Golomb67}
\label{state-diagram1-7}
Let $F(x_{1},x_{2},\ldots, x_{n})$ be the feedback function of an $n$-stage FSR.
If its state diagram only consists of disjoint cycles, then this FSR is a nonsingular FSR.
Any two cycles $C_i$ and $C_j$ can be joined into a single cycle when the successors of $Y^n\in C_i$ and its conjugate $\widehat{Y^n}\in C_j$ are interchanged, and a de Bruijn cycle   be obtained finally.
\end{lemma}

\begin{lemma}
\cite{Golomb67}
\label{state-diagram1-8}
 Let $F(x_{1},x_{2},\ldots, x_{n})$ be the feedback function of an $n$-stage nonsingular FSR.
 If any state $Y^n$ and its conjugate $\widehat{Y^n}$ belong to the same cycle, then the state diagram of this FSR is a de Bruijn cycle.
\end{lemma}

\begin{algorithm}[h!]
\caption{Find  conjugate pairs $conj(Z^n,\widehat{Z^n})$ between any two adjacent cycles in $C$}
\label{alg1-2}
\begin{algorithmic}[1]
\REQUIRE  $P=P_1\cup P_1\cup\ldots\cup P_t$ and $C=\{C_1,C_2,\ldots,C_t\}$ determined by Algorithm  \ref{alg1-1}.
\ENSURE The  conjugate pairs $conj(Z^n,\widehat{Z^n})$ between any two distinct cycles  $C_i,C_j\in C$.
\STATE For $P_i=\{p_{i,1},p_{i,2},\ldots,p_{i,s_i}\}$ and  $P_j=\{p_{j,1},p_{j,2},\ldots,p_{j,s_j}\}$,
\IF {$p_{i,d}$ and $p_{j,e}$ are two $(n-2)$-length paths and $Y^n_{i,d},Y^n_{j,e}\not\in S^n_{100}$,}
\STATE there is only one conjugate pair $conj(Z^n,\widehat{Z^n})$ between  $p_{i,d}$ and $p_{j,e}$, where $Y^n_{i,d}=t(p_{i,d}),Y^n_{j,e}=t(p_{j,e})$.
\ENDIF
\STATE  For any $p_{i,w}\in P_i$ and  $p_{j,z}\in P_j$,
\IF { $Y^n_{i,w}$ and $Y^n_{j,z}$ belong to the same set of $S^n_{100}$, $S^n_{001}$, $S^n_{010}$ and $S^n_{111}$,}
\STATE there is at most one conjugate pair $conj(Z^n,\widehat{Z^n})$ between  $p_{i,w}$ and  $p_{j,z}$, where $Y^n_{i,w}=t(p_{i,w}),Y^n_{j,z}=t(p_{j,z})$.
\ELSE
\RETURN $C_i$ and $C_j$ are not adjacent.
\ENDIF
\RETURN {$C_F=\{conj(Z^n,\widehat{Z^n})\}$, where $Z^n$ and $\widehat{Z^n}$ are the conjugate states between two different cycles, respectively.}
\end{algorithmic}
\end{algorithm}

According to Lemma \ref{state-diagram1-7} and \ref{state-diagram1-8}, the conjugate pairs between any two cycles  are  needed in the construction of de Bruijn sequences by the cycle joining method.
The method  to find conjugate pairs $conj(Z^n,\widehat{Z^n})$ between  two any cycles in $C$ is given in  Algorithm \ref{alg1-2}.
Since each  trigeminal vertex has two conjugate predecessor states, the conjugate pairs $conj(Z^n,\widehat{Z^n})$ appear at one trigeminal vertex nearby.
Through the characteristics of state diagram, we notice that a  conjugate pair consists of either one state of one cycle and  one state adjacent to the cycle, or two  states of  the same  perfect  binary directed tree.

Let $P=P_1\cup P_1\cup\ldots\cup P_t$ and $C=\{C_1,C_2,\ldots,C_t\}$ be the set of directed paths and the set of cycles obtained by Algorithm \ref{alg1-1}, respectively, where $C_{i}: [p_{i,1}\rightarrow p_{i,2}\rightarrow\ldots\rightarrow p_{i,s_{i}}\rightarrow p_{i,1}]$ and $P_{i}=\{p_{i,1},p_{i,2},\ldots,p_{i,s_{i}}\}(1\leq i\leq t)$.

\begin{theorem}
\label{state-diagram1-9}
Given a  directed path searching in the state diagram as in Algorithm \ref{alg1-1} with feedback function $F(x_{1},x_{2},\ldots, x_{n})=x_{n-1}+x_n \ (n\geq3)$, if we repeat the searching conjugate pair $conj(Z^n,\widehat{Z^n})$ between two any cycles in $C$ as in Algorithm \ref{alg1-2},  then  conjugate pairs $conj(Z^n,\widehat{Z^n})$ between  two any cycles will be obtained.
\begin{proof}
Since the cycle obtained by the Algorithm \ref{alg1-1} is a union of  some directed paths,  the problem of considering the conjugate pair $conj(Z^n,\widehat{Z^n})$ between the two cycles is  equivalent to the problem of considering the conjugate pair $conj(Z^n,\widehat{Z^n})$ between their directed paths.
Note that a  conjugate pair $conj(Z^n,\widehat{Z^n})$ consists of either one state of one cycle and  one state adjacent to the cycle, or two  states of  the same  perfect  binary directed tree.

Case 1: When the last states of two directed paths  belong to the same cycle of the original state diagram, that is, when two $(n-2)$-length directed paths have their last states  in the same cycle, they can only share one conjugate pair $conj(Z^n,\widehat{Z^n})$.
In fact, the last state of one path and the last but one state of the other path are a pair of conjugate states.
Since the state diagram only contains two cycles  $[0]$ and $[0,1,1]$, we only consider the $(n-2)$-length paths whose tails are in the set $S_0^n\setminus\{\mathbf{0}^n\}$.

Case 2: There are only  conjugate pairs $conj(Z^n,\widehat{Z^n})$ between the two directed paths  of the two cycles whose tails belong to the same set of $S_{100}^n$, $S_{001}^n$, $S_{010}^n$ and $S_{111}^n$.

Based on the set of directed paths and the set of cycles obtained by Algorithm \ref{alg1-1},  conjugate pairs $conj(Z^n,\widehat{Z^n})$ between  two any cycles  be obtained.
\end{proof}
\end{theorem}

The number of directed paths obtained by Algorithm \ref{alg1-1} is $2^{n-1}$.
The Algorithm \ref{alg1-2} is to find conjugate pairs $conj(Z^n,\widehat{Z^n})$ by comparing the directed paths between different cycles.
The time complexity of Algorithm \ref{alg1-2} is $O(2^{n-1})$.
For the set of cycles obtained by Algorithm \ref{alg1-1} and conjugate pairs $conj(Z^n,\widehat{Z^n})$ obtained by Algorithm \ref{alg1-2}, a new class of de Bruijn cycles  be obtained with the  cycle-joining method.

\section{Example}
\label{sec1-5}
\begin{figure}[h!]
\vspace{-0.5cm}
\centering
\subfigure[$G_1$]{
\label{Fig.sub.1-1}
\includegraphics[width=8.00cm,height=7.00cm]{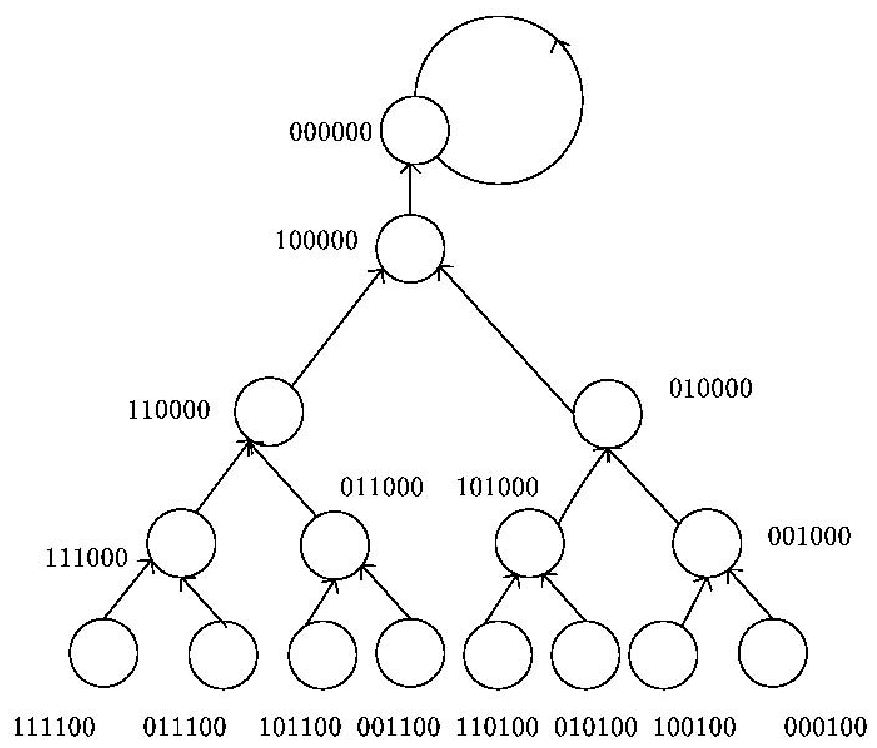}}
\subfigure[$G_2$]{
\label{Fig.sub.1-2}
\includegraphics[width=7.00cm,height=7.00cm]{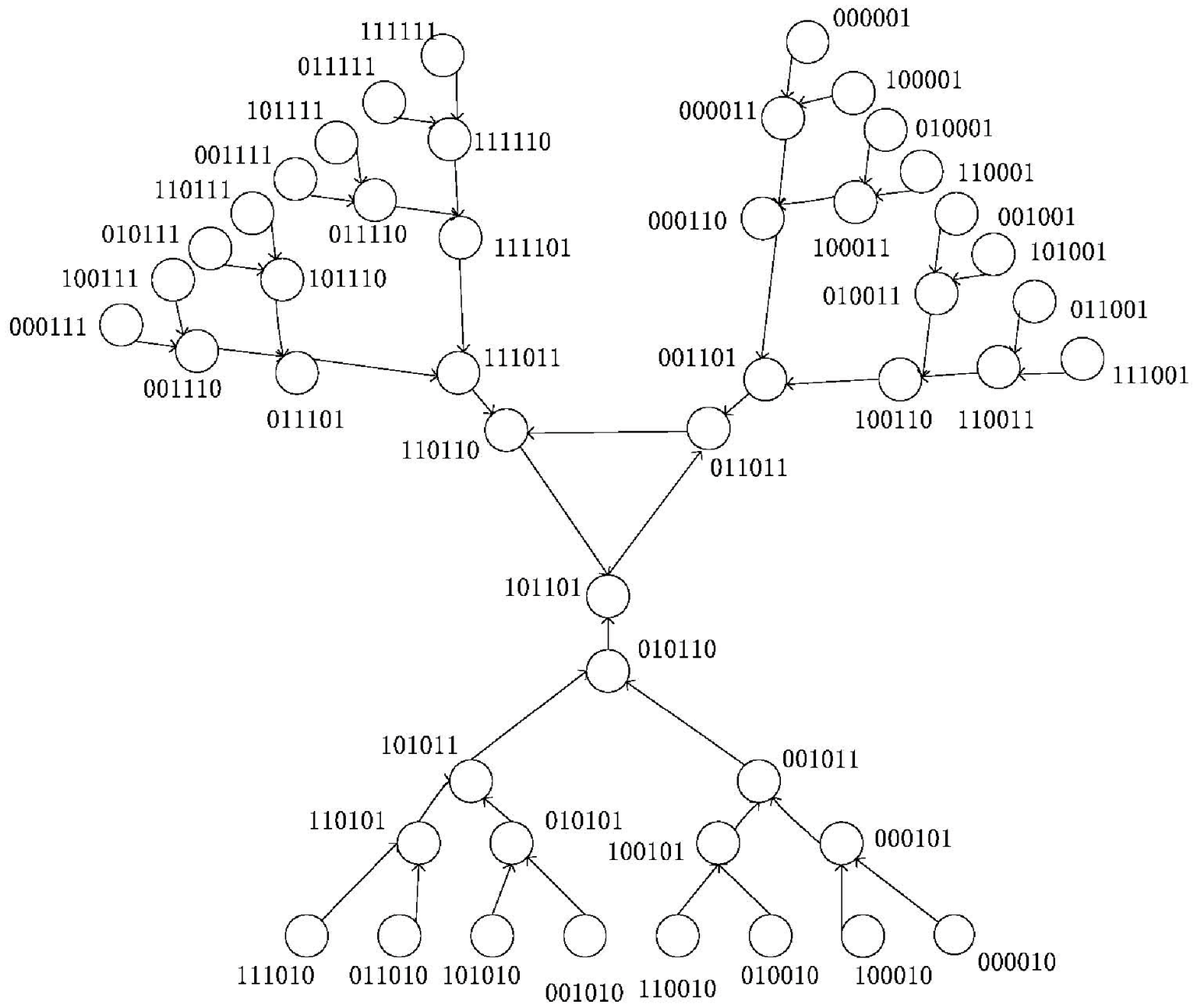}}
\vspace{-0.5cm} \caption{The state diagrams  of   $F(x_{1},x_{2},\ldots, x_6)=x_5+x_6 $}
\label{Fig.lable1-1}
\end{figure}

Let $F(x_1,x_2,\ldots,x_6)=x_5+x_6$, the state diagram  of   $F(x_{1},x_{2},\ldots, x_6)=x_5+x_6 $ is given in Figure~\ref{Fig.lable1-1}.

\begin{table*}[!b]
\vspace{-0.5cm}
\caption{Repeat the directed path searching in the state diagram as in Algorithm \ref{alg1-1}}
\label{tab1-1}
\begin{center}
\begin{tabular}{|c|c|c|}\hline
$Y_{1,1}^6=(1,1,1,1,1,1)$&$l_{1,1}=4$&$p_{1,1}:[Y_{1,1}^6\rightarrow LY_{1,1}^6\rightarrow\ldots\rightarrow L^4Y_{1,1}^6]$\\\hline
$Y_{1,2}^6=(1,0,1,1,0,0)$&$l_{1,2}=4$&$p_{1,2}:[Y_{1,2}^6\rightarrow LY_{1,2}^6\rightarrow\ldots\rightarrow L^4Y_{1,2}^6]$\\\hline
$Y_{1,3}^6=(0,0,0,0,0,1)$&$l_{1,3}=4$&$p_{1,3}:[Y_{1,3}^6\rightarrow LY_{1,3}^6\rightarrow\ldots\rightarrow L^4Y_{1,3}^6]$ \\\hline
$Y_{1,4}^6=(1,1,0,1,1,1)$&$l_{1,4}=3$&$p_{1,4}:[Y_{1,4}^6\rightarrow LY_{1,4}^6\rightarrow L^2Y_{1,4}^6]$\\\hline
$Y_{1,5}^6=(1,1,1,0,1,0)$&$l_{1,5}=4$&$p_{1,5}:[Y_{1,5}^6\rightarrow LY_{1,5}^6\rightarrow\ldots\rightarrow L^4Y_{1,5}^6]$\\\hline
$Y_{1,6}^6=(0,1,1,0,1,0)$&$l_{1,6}=1$&$p_{1,6}:[Y_{1,6}^6]$\\\hline
$Y_{1,7}^6=(1,1,0,1,0,0)$&$l_{1,7}=3$&$p_{1,7}:[Y_{1,7}^6\rightarrow LY_{1,7}^6\rightarrow L^2Y_{1,7}^6]$\\\hline
$Y_{1,8}^6=(1,0,0,0,0,1)$&$l_{1,8}=1$&$p_{1,8}:[Y_{1,8}^6]$\\\hline
$Y_{1,9}^6=(0,0,0,0,1,0)$&$l_{1,9}=3$&$p_{1,1}:[Y_{1,9}^6\rightarrow LY_{1,9}^6\rightarrow L^2Y_{1,9}^6]$\\\hline
$Y_{1,10}^6=(0,1,0,1,1,1)$&$l_{1,10}=1$&$p_{1,2}:[Y_{1,10}^6]$\\\hline
$Y_{1,11}^6=(1,0,1,1,1,1)$&$l_{1,11}=2$&$p_{1,3}:[Y_{1,11}^6\rightarrow LY_{1,11}^6]$ \\\hline
$Y_{1,12}^6=(1,1,1,1,0,0)$&$l_{1,12}=2$&$p_{1,4}:[Y_{1,12}^6\rightarrow LY_{1,12}^6]$\\\hline
$Y_{1,13}^6=(1,1,0,0,0,1)$&$l_{1,13}=2$&$p_{1,13}:[Y_{1,13}^6\rightarrow LY_{1,13}^6]$\\\hline
$Y_{1,14}^6=(0,0,0,1,1,1)$&$l_{1,14}=2$&$p_{1,14}:[Y_{1,14}^6\rightarrow LY_{1,14}^6]$\\\hline
$Y_{1,15}^6=(0,1,1,1,0,0)$&$l_{1,15}=1$&$p_{1,15}:[Y_{1,15}^6]$\\\hline
$Y_{1,16}^6=(1,1,1,0,0,1)$&$l_{1,16}=3$&$p_{1,16}:[Y_{1,16}^6\rightarrow LY_{1,16}^6\rightarrow L^2Y_{1,16}^6]$\\\hline
$Y_{1,17}^6=(0,0,1,1,0,0)$&$l_{1,17}=1$&$p_{1,17}:[Y_{1,17}^6]$\\\hline
$Y_{1,18}^6=(0,1,1,0,0,1)$&$l_{1,18}=1$&$p_{1,18}:[Y_{1,18}^6]$\\\hline
$Y_{1,19}^6=(1,1,0,0,1,0)$&$l_{1,19}=2$&$p_{1,19}:[Y_{1,19}^6\rightarrow LY_{1,19}^6]$ \\\hline
$Y_{1,20}^6=(0,0,1,0,1,0)$&$l_{1,20}=2$&$p_{1,20}:[Y_{1,20}^6\rightarrow LY_{1,20}^6]$\\\hline
$Y_{1,21}^6=(1,0,1,0,1,0)$&$l_{1,21}=1$&$p_{1,21}:[Y_{1,21}^6]$\\\hline
$Y_{1,22}^6=(0,1,0,1,0,0)$&$l_{1,22}=1$&$p_{1,22}:[Y_{1,22}^6]$\\\hline
$Y_{1,23}^6=(1,0,1,0,0,1)$&$l_{1,23}=2$&$p_{1,23}:[Y_{1,23}^6\rightarrow LY_{1,23}^6]$\\\hline
$Y_{1,24}^6=(1,0,0,1,1,1)$&$l_{1,24}=1$&$p_{1,24}:[Y_{1,24}^6]$\\\hline
$Y_{1,25}^6=(0,0,1,1,1,1)$&$l_{1,25}=1$&$p_{1,25}:[Y_{1,25}^6]$\\\hline
$Y_{1,26}^6=(0,1,1,1,1,1)$&$l_{1,26}=1$&$p_{1,26}:[Y_{1,26}^6]$\\\hline\hline
$Y_{2,1}^6=(1,0,0,1,0,0)$&$l_{2,1}=2$&$p_{2,1}:[Y_{2,1}^6\rightarrow LY_{2,1}^6]$ \\\hline
$Y_{2,2}^6=(0,1,0,0,0,1)$&$l_{2,2}=1$&$p_{2,2}:[Y_{2,2}^6]$\\\hline
$Y_{2,3}^6=(1,0,0,0,1,0)$&$l_{2,3}=1$&$p_{2,3}:[Y_{2,3}^6]$\\\hline
$Y_{2,4}^6=(0,0,0,1,0,0)$&$l_{2,4}=1$&$p_{2,4}:[Y_{2,4}^6]$\\\hline
$Y_{2,5}^6=(0,0,1,0,0,1)$&$l_{2,5}=1$&$p_{2,5}:[Y_{2,5}^6]$\\\hline
$Y_{2,6}^6=(0,1,0,0,1,0)$&$l_{2,6}=1$&$p_{2,6}:[Y_{2,6}^6]$\\\hline\hline
$C_{1}$&&$ [p_{1,1}\rightarrow p_{1,2}\rightarrow\ldots\rightarrow p_{1,26}\rightarrow p_{1,1}]$\\\hline
$C_{2}$&&$ [p_{2,1}\rightarrow p_{2,2}\rightarrow\ldots\rightarrow p_{2,6}\rightarrow p_{2,1}]$\\\hline
\end{tabular}
\end{center}
\vspace{-0.5cm}
\end{table*}

Step 1:  Start a search at one of the leaves $Y_{1,1}^6=(1,1,1,1,1,1)$, and trace a directed path downward in the tree.
Record the path and all vertices of the path from the tree.
Then repeat the directed path search from the state diagram until the forest is empty.
Note that each isolated vertex is considered as a path of length $0$, and we finally get $32$ different directed paths and $2$ cycles.
In Table \ref{tab1-1}, we give the searching way.

Step 2:  In Table \ref{tab1-2}, the classification of the tails of paths of Algorithm \ref{alg1-1} is given.

Step 3:  In Table \ref{tab1-3}, $5$ conjugate pairs $conj(Z^n,\widehat{Z^n})$  between two cycles  $C_1$ and $C_2$ are found.
Then interchange their successors, and $5$ de Bruijn cycles can be obtained.

\begin{table*}[!t]
\vspace{-0.5cm}
\begin{center}
\caption{ The classification of the tails of paths}
\label{tab1-2}
\begin{tabular}{|c|c|c|}
\hline
 & $C_1$&$C_2$ \\\hline
 $S_{111}^6$&$Y_{1,1}^6,Y_{1,4}^6,Y_{1,10}^6,Y_{1,11}^6,Y_{1,14}^6,Y_{1,24}^6,Y_{1,25}^6,Y_{1,26}^6$& -----\\\hline
 $S_{100}^6$&$Y_{1,2}^6,Y_{1,7}^6,Y_{1,12}^6,Y_{1,15}^6,Y_{1,17}^6,Y_{1,22}^6$& $Y_{2,3}^6,Y_{2,6}^6$\\\hline
 $S_{001}^6$&$Y_{1,3}^6,Y_{1,8}^6,Y_{1,13}^6,Y_{1,16}^6,Y_{1,18}^6,Y_{1,23}^6$& $Y_{2,2}^6,Y_{2,5}^6$\\\hline
 $S_{010}^6$&$Y_{1,5}^6,Y_{1,6}^6,Y_{1,9}^6,Y_{1,19}^6,Y_{1,20}^6,Y_{1,21}^6$&$Y_{2,1}^6,Y_{2,4}^6$ \\\hline
\end{tabular}
\end{center}
\vspace{-0.5cm}
\end{table*}

\begin{table*}[!t]
\vspace{-0.5cm}
\begin{center}
\caption{Conjugate pairs   between two cycles  $C_1$ and $C_2$}
\label{tab1-3}
\begin{tabular}{|c|c|}
 \hline
     \multicolumn{2}{| c |}{$conj(Z^n,\widehat{Z^n})$ }\\\hline
 $(1,0,1,0,0,0)$& $(0,0,1,0,0,0)$ \\\hline
 $(1,1,0,0,0,1)$& $(0,1,0,0,0,1)$\\\hline
 $(0,0,0,0,1,0)$& $(1,0,0,0,1,0)$\\\hline
 $(1,0,1,0,0,1)$& $(0,0,1,0,0,1)$\\\hline
 $(1,1,0,0,1,0)$& $(0,1,0,0,1,0)$ \\\hline
\end{tabular}
\end{center}
\vspace{-0.5cm}
\end{table*}

\section{Conclusion}
\label{sec1-6}
The state diagrams  of a class of singular LFSR are discussed.
Some properties of these  singular linear feedback shift registers are also given.
An  algorithm is presented to construct a new class of de Bruijn cycles  from the state diagrams of these singular LFSR.
This is the first time to construct  de Bruijn cycles based on  singular linear shift registers.
In this method,  cycle  structures are obtained by modifying the state diagrams firstly, then the conjugate pairs between  cycles are searched in the directed paths set.
The de Bruijn cycles are realized by using the cycle-joining method finally.

\section*{Acknowledgements}

This work is supported by the National Natural Science Foundation of China under grants 61672414, and the
National Cryptography Development Fund under Grant MMJJ20170113.

\end{document}